\documentclass[preprint]{article}

\usepackage{array}

\usepackage{float}
\usepackage{placeins}

\usepackage{subfigure}

\usepackage{amsmath, amssymb, graphicx}

\def \R {\mathbb R}
\def \s {\sigma}
\def \S {\Sigma}

\def \beg{\begin{eqnarray}}
\def \en{\end{eqnarray}}
\def \be*{\begin{eqnarray*}}
\def\e*{\end{eqnarray*}}
\def \di{\displaystyle}
\def\bit{\begin{itemize}}
\def \eit{\end{itemize}}
\def \E{\mathbb E}
\def \N{\mathbb N}

\def \w{\widehat}

\def \beg{\begin{eqnarray}}
\def \en{\end{eqnarray}}
\def \be*{\begin{eqnarray*}}
\def \e*{\end{eqnarray*}}
\def \di{\displaystyle}
\def \bit{\begin{itemize}}
\def \eit{\end{itemize}}

\renewcommand{\geq}{\geqslant}
\renewcommand{\hat}{\widehat}

\renewcommand{\P}{\mathbb{P}}
\def\argmax{\mathop{\mbox{\sl\em argmax}}}

 \def \cqfd{\hspace*{14.6cm}$\blacksquare$}

\newtheorem{thm}{Theorem}

\newtheorem{rmk}{Remark}
\newtheorem{prop}{Proposition}

\title{Comparing two contaminated samples}

\author{D. Pommeret (pommeret@univmed.fr) 
\\ Institute of Mathematics of Luminy - Case 907 \\
Aix-Marseille II University - 13288 Marseille cedex 9,  France}

\date{}
\begin{document}

\maketitle




\begin{abstract}
In this paper we consider
the  problem of testing whether two samples of contaminated data, possibly paired,  are from the same distribution.
Is is assumed that the contaminations are additive noises with known moments of all orders.
The test statistic is based on the polynomials moments
of the difference between observations and  noises.  .
A data driven selection  is proposed to choose automatically the number of involved polynomials.
We present a simulation study in order to investigate the power of the proposed test within discrete and continuous cases.
A real-data example is presented to demonstrate the method.
\end{abstract}

\textbf{keyword}
 contaminated data; data-driven; 
two sample test

\section{Introduction}





The classical two-sample problem concerning i.i.d. observations has been extensively  studied in the literature.
We propose in this paper to extend this problem to the case of two contaminated samples when a noise is added to each sample.
More precisely, we consider two samples, $X_1, \cdots, X_n$ and $U_1, \cdots, U_k$, from the following two models
\beg
\label{convol}
X =Y+Z, & {\rm and } &  U= V+W,
\en
where $Y$ and $Z$ (resp. $V$ and $W$) are two independent
random variables.
 It is also assumed that $Z$ and $W$ are independent. However this paper concerns independent as well as paired variables $X$ and $U$ since $Y$ and $V$ can be dependent.  We keep this hypothesis
through the paper putting $n=k$
(the more general  case being  easily obtained).
We  assume that all moments of $Z$ and $W$ exist and are known.
We are interested in testing the equality of the distribution of $Y$ and $V$.
Our aim  is to construct an omnibus test for the general non parametric hypothesis
\beg
\label{hyp}
H_0: {\cal L}_Y= {\cal L}_V & {\rm against} & H_1: {\cal L}_Y\neq {\cal L}_V,
\en
where ${\cal L}_Y$ and ${\cal L}_V$ refers to the distribution of $Y$ and $V$.
For that we extend the one-sample smooth test inspired of Neyman (1937) (see also Rayner and Best, 1989, for a general introduction) to the two-sample case under (\ref{convol}).
 For the one sample problem, the smooth test is an  omnibus approach which consists in coming down to  parametric hypotheses. Then the smooth statistic is composed of different elements each able to detect a departure from the null hypothesis. This approach can  be naturally   extended to the two sample case, as in Rayner and Best (2001)
 (see also Chervoneva and Iglewicz, 2005). In addition, Ledwina (1994) introduced a data driven procedure permitting
 to select automatically the number of elements of the statistic. The automatic selection is based on the Schwarz (1978) criterion. Janik-Wr\'oblezca and Ledwina (2000) first used this technique combined with rank statistic for the two sample problem.
  Recently   Ghattas { et al.} (2011) obtained a data driven  test for the two paired sample problem.
Various extensions of the data driven smooth test have been proposed, particularly in the context of survival data in Krauss (2009) when samples are right censored, or in the context of  detection of changes in  Antoch { et al.} (2008)  reducing the problem to a two sample subproblem.

From (\ref{convol}) it is clear that the unknown moments of $Y$ (resp. $V$) can be expressed
in terms of moments of $X$ and $Z$ (resp. $U$ and $W$).
 The proposed smooth test  is  based on the difference between the $k$ first moments of $Y$ and $V$.
The order $k$ determines the number  of components of the test statistic.
We then  adapt the data driven approach
permitting to select automatically this number. We first consider the case where $k$
varies between $1$ and $K$, for $K$ a fixed integer. Then
 we let $k$ tend to infinity more slowly than the sample size.
For asymptotic results we make  an assumption  on the smallest eigenvalue of the sample covariance matrix.
But in practice, the data driven procedure is effective in the first case with $K$  fixed large enough, as shown in our simulations.
Finally, we apply our method to the UEFA champion's league data from Meintanis  (2007).

Before describing our test procedure  we offer a few examples that illustrate the situation (\ref{convol}).
\bit
\item
\emph{Evaluation by experts}.  During an assessment, such as sensory analysis, it is very common that experts are biased in their judgments. This bias is commonly observed and assessed during training and can be assumed to be known in distribution. Typically, one can assume a normal distribution with mean and variance associated with each expert. In this case, if we want to compare the distribution of two products evaluated by two experts, we are reduced to the situation (\ref{convol}) where $X$ and $U$ coincide with the two experts scoring, $Z$ and $W$ being their errors.

\item
\emph{Ruin theory.}
Another situation that can be encountered in ruin theory  is  the random sum of claims, $\sum e_i$, where $e_i$ are i.i.d. random variables with known exponential distribution. The number of claims can be decomposed into a fixed known value,  $n$, and a random value,   $N$, representing an aggregation of different claims. Thus if we observe  two sums
\be*
X  =
\di\sum_{i=1}^{n_1+N_1}e_i
& {\rm and } &
U  =
\di\sum_{i=1}^{n_2+N_2}v_i,
\e*
where $e_i$ and $v_i$ are i.i.d., one problem is to compare the randomness structure $N_1$ and $N_2$, that is to test the equality of the distributions of these two  variables. This problem coincides with (\ref{convol}) since it is equivalent to testing the equality of the distributions  of
$Y=\sum_{i=1}^{N_1}e_i$ and $V= \sum_{i=1}^{N_2}e_i$.
\item
\emph{Mixture model}.
The deconvolution problem  is also
related to a mixture problem since a particular case of (\ref{convol}) is  the
location mixture situation of the form
\be*
f_X(x)=\int f_Y(x-m) f_Z(dm), & &
f_U(x)=\int f_V(x-m) f_W(dm),
\e*
 with $m$ the location
parameter, $f_Y$, $f_V$ the unknown mixed densities and $f_Z$, $f_W$ the known mixing densities. This situation can be encountered when  finite mixture distributions have known components, and when the purpose is to compare their associated sub-populations associated with these components. We can also reverse the roles of $Y,V$ and $Z,W$ and be interested in the comparison of two linear mixed models with Gaussian noise and unknown random effects.
 \item
\emph{Extreme values}.
Contaminated model  can be also viewed as a model for extremal values considering the convolutions
\be*
X=\alpha Y+Z, & &
U=\beta V+W,
\e*
 where $\alpha$ and $\beta$ are Bernoulli with small parameter representing the occurrence of an extreme event.  Often the non-extreme distributions  of $Z$ and $W$ are well observed and known and we can be interested in the comparison of the extreme distributions of $Y$ and $V$. Assuming that one knows when these rare events occur, they are observed with a known noise as in (\ref{convol}).
\item
\emph{Scale mixture}.
Finally, it is current to observe the product of two variables, say
\be*
X=YZ , & & U=VW.
\e*
For instance, that is the case for Zero Inflated distributions, when $Y$ and $V$ are Bernoulli random variables and $Z$ and $W$ are discrete random variables. Without loss of generality, by translating all variables, we can use a log-transformation to recover (\ref{convol}). Many other cases can be envisaged as $X=Z/Y$ and $U=W/V$ with $Z$ and $W$ normally distributed.
\eit

The  paper is organized as follows. In Section 2 we introduce the method based on
polynomial expansions for testing the equality of the
two contaminated densities. In Section 3  we
propose a simple data driven procedure  that we extend  to the case where the number of components of the statistic tends to infinity, with additional assumptions. In Section 4, finite-sample properties of the proposed test
statistics are examined through Monte Carlo simulations. The analysis of a the champion's league data set is provided in
Section 5. Section 6  contains a brief discussion.

\section{Statistical method}

Consider simultaneously two (possibly paired) samples $X_1,\cdots, X_n$ and $U_1,\cdots, U_n$ following (\ref{convol})
and such that all moments exist and characterize the associated  distributions.
 Is is assumed that
 the moments of $Z$ and $W$ are known. From (\ref{convol})  we have the following  two expansions for all integer $i$
\beg
\E(X^i) =  \di\sum_{j=0}^i
c_{ij}\E(Y^j)z_{i-j},
& {\rm and } &
\E(U^i)  =  \di\sum_{j=0}^i
c_{ij}E(V^j)w_{i-j},
\label{basic}
\en
with $c_{ij}=(i!)/((j!)((i-j)!))$, $z_{i-j}=\E(Z^{i-j})$, and  $w_{i-j}=\E(W^{i-j})$.
Write $a_i=\E(Y^i)$ and $b_i=\E(V^i)$. The null hypothesis coincides with $a_i=b_i$, $\forall i=1,2,\cdots$, and our testing procedure reduces to the parametric testing problem: $\forall i=1,\cdots, k, $ $a_i-b_i=0$, when  $k$ gets large. We shall let $k$ tend to infinity, with a speed depending of the sample size, and its choice will be done automatically  by a data driven method.
Inverting (\ref{basic}) we get
\beg
a_i=  \E(P_i(X))
& {\rm and } &
b_i  =  \E(Q_i(U)),
\label{basic2}
\en
where $P_i$ and $Q_i$ are polynomials of degree $i$.
For instance the first three terms are
\be*
P_1(x) & = & x-z_1,
\\
P_2(x) & = & x^2-2z_1P_1(x)-z_2,
\\
P_3(x) & = & x^3 -3z_1P_2(x)-3z_2P_1(x)-z_3.
\e*
To construct the test statistics we consider the vector of differences
\be*
V_s(k)&=& \left(P_i(X_s) -Q_i(U_s)\right)_{1 \leq i \leq k},
\e*
and we put
\be*
J_n(k)&=&\frac{1}{\sqrt{n}}\sum_{s=1}^nV_s(k).
\e*
 Under $H_0$, $J_n(k)$ has mean zero and finite  $k\times k$ variance-covariance matrix
 $$\S(k)=\E_0\big (V_1(k) V_1(k)^{'}\big),$$ where $\E_0$ denotes
 the expectation under $H_0$ and $V_1(k)^{'}$ is the transposition of $V_1(k)$.
  Next, let  us define the empirical version of $\S(k)$ under $H_0$, that is the $k\times k$ matrix
  $$\hat{\S}_n(k)=\frac{1}{n}\sum_{s=1}^nV_s(k) V_s(k)^{'}.$$   In the
  following, we assume that  $\hat \S_n(k)$ is a positive-definite matrix so that the
  corresponding inverse matrix and its square root exist. Note that this condition is satisfied a.s. for $n$ large enough since the estimator is consistent.
 We consider  the test statistic
\begin{equation} \label{Tk}
T_n(k)=J_n(k)^{'}\hat{\S}_n(k)^{-1}J_n(k)=\|\hat{\S}_n(k)^{-1/2} J_n(k)\|^2,
\end{equation}
 where
$\|. \|$  denotes the euclidian norm on
$\mathbb R^k$.
Application of the Central Limit Theorem shows that under $H_0$, $T_n(k)$ converges in distribution to a
$\chi^2$  random variable with $k$ degrees of freedom as $n$ tends to infinity.
The strategy is to select an appropriate degree $k$; that is, a correct number of components in the
test statistics.
In addition, observe that the null hypothesis can be rewritten
as   $H_0: \theta=0$
where $\theta=\E(V_1(k))$. Suppose that the maximum likelihood estimator $\widehat \theta$ of $\theta$ equals the
empirical mean of the sample of the $V_s(k)$'s, that is  $\widehat \theta=\frac{J_n(k)}{\sqrt{n}}$, as it is the case for instance when
the distribution of $V_1(k)$ belongs to an exponential family. Then, $T_n(k)$ is the score  statistic and the Schwarz criteria is well adapted to get an automatic selection of $k$.

\section{Data driven approach}

In  this section, the data-driven method introduced by Ledwina (1994) (see also
Inglot et al.  1997) is used to optimize the parameter $k$ in our test statistic. It is based on a
modified version of Schwarz's Bayesian information rule.
The optimal value of $k$, denoted by $S_n$, is such that
\begin{equation} \label{rule}
S_n= \min\big \{\argmax_{1 \leq k \leq d(n)} ( T_n(k) - k\log(n))  \big\},
\end{equation}
where $d(n)$ can be either  fixed, equal to $K$, or increasing such that  $\lim_{n\rightarrow\infty}d(n)=\infty$.
Once $S_n$ is determined, the test statistic is applied with $k=S_n$. More precisely, we
 use for our testing problem   the  statistic $T_n(S_n)$.
Hereafter,   the asymptotic distribution of the test statistic is derived  under
the null hypothesis for cases where $d(n)$ is fixed or unbounded.

\subsection{The case where $d(n)=K$ is fixed.}
\begin{thm}
\label{theo0}
Assume that  $d(n) = K\geq 1$ is fixed. Under $H_0$, when $n$ tends to infinity,
$T_n(S_n)$ converges in distribution to a
$\chi^2$ random variable with 1 degree of freedom.
\end{thm}
The
proof is fairly standard and follows Ledwina (1994). We will detail a more general  proof
 in the case where $d(n)$ is unbounded (see Theorem \ref{theo2}).
\begin{rmk}
In our simulations, we fixed $K$ large enough, in the sense that its value was neither reached by $S_n$, either under the null (for empirical level calculations) or under alternatives (for empirical power calculations).
\end{rmk}
\subsection{The case where $d(n)$ is unbounded.}
 Let us denote
by $\P_0$ and $\E_0$ the probability and the expectation
under $H_0$.
 Write $\widehat \lambda_{\min}(k)$
the smallest
eigenvalue of $\hat \S_n(k)$.
We now let $d(n)$ tend to infinity under  the following two conditions:
\begin{description}
\item{{(A1)}}
 $d(n)^2/ \E_0\big({\w{\lambda}_{\min}(d(n))}\big) = o_{\P_0}(\log(n))$.
 \item{{(A2)}} There exists some positive constant $M$ such that for all $k>0$,
 $$
 \di\frac{1}{k}\di\sum_{i=1}^k \E_0\bigg( Z_i^4 \bigg) <M.$$
 where $Z_i = P_i(X)-Q_i(U)$.
\end{description}

\begin{rmk}
The condition (A1) can be compared to results obtained in the framework of random matrices.
For instance, Bai and Yin (1993) (see also Silverstein, 1985, for the particular Gaussian case) considered  the case where the entries
$Z_{ij}=P_i(X_j)-Q_i(X_j)$ are independent and identically
distributed with finite fourth moment (this moment condition may be compared with (A2)). They shown that almost surely
$\lim \w{\lambda}_{\min}(d(n))=1$ when $d(n)/n \rightarrow 0$. Then when the random series $\w{\lambda}_{\min}(d(n))$ is bounded  we get $\lim \E\big(\w{\lambda}_{\min}(d(n))\big)=1$ and $d(n)$ can be chosen as $o_{\P_0}(\sqrt{\log(n)})$.

Assumption (A2) states that the fourth moment is bounded on average.
It is similar to Assumption 2 stated  in  Ledoit and Wolf (2004).
More precisely, Ledoit and Wolf used a condition on the eighth moment which is somewhat  more restrictive.

\end{rmk}

\begin{thm}
\label{theo2}
Let assumptions (A1) and (A2) hold.
 Thus, under $H_0$,
$T_n(S_n)$ converges in distribution to a
$\chi^2$ random variable with 1 degree of freedom.
\end{thm}
\textbf{Proof}
The proof is partly inspired by Janic-Wr\'oblewska and  Ledwina (2000).
First note that the greatest eigenvalue
of $\hat \S_n(k)^{-1}$  is the inverse of its smallest eigenvalue.
Then we have $\||\hat \S_n(k)^{-1}\|| = 1/\hat \lambda_{\min}(k)$, where $\|| . \||$ stands for the spectral norm.
 Under $H_0$, it is clear that   $T_n(1)$ converges to a $\chi^2$ random
 variable
with one degree of freedom. Then we have  to prove that
 $\P_0(S_n = 1)$ tends to 1
as $n$ tends to infinity, or equivalently that  $\P_0(S_n \geq 2)$ tends to 0. Let us set $a_n(k)=(k-1)\log{n}$.
By definition of $S_n$, we have
\begin{equation}\label{e1}
\P_0(S_n\geq 2)=\sum_{k=2}^{d(n)}\P_0(S_n=k) \leq
\sum_{k=2}^{d(n)}\P_0\bigg(T_n(k)^{1/2} \geq
\sqrt{a_n(k)}\bigg).
\end{equation}
Using the standard norm's inequalities for matrices and vectors we get
 \begin{equation} \label{mat4}
 T_n(k)= { J_n(k)'\hat \S_n(k)^{-1}J_n(k)}\leq  {\|J_n(k)\|^2}\||\hat \S_n(k)^{-1}\||,
 \end{equation}
that we combine with Markov inequality to obtain
\begin{eqnarray*} \label{ee3}
\P_0\bigg (T_n(k)^{1/2}\geq {\sqrt{a_n(k)}}\bigg)
&\leq& \P_0\bigg ({\|J_n(k)\|}\||\hat \S_n(k)^{-1}\||^{1/2} \geq \sqrt{{a_n(k)}}
\bigg )
\\
& \leq &
 {\frac{\E_0 \big (\| J_n(k)\| \||\S_n(k)^{-1}\||^{1/2}\big )}{\sqrt{a_n(k)}}}
\\
& \leq &
 \di\frac{\bigg(\E_0 \big (\| J_n(k)\|^2\big) \E_0\big(\||\S_n(k)^{-1}\||\big)\bigg)^{1/2}}{\sqrt{a_n(k)}}
\\
& = &
\di\frac{\bigg(\E_0 \big (\| J_n(k)\|^2\big)\bigg)^{1/2}}{ \E_0\big(\w{\lambda}_{\min}(k)\big)^{1/2}\sqrt{a_n(k)}}.
 \end{eqnarray*}
 Using the independence of the pairs $(X_s,Y_s)_{1
\leq s\leq n}$,  we get
 \beg
  \E_0 \bigg (\|
J_n(k)\|^2\bigg )&=&\E_0 \bigg (\frac{1}{n}
\sum_{s=1}^n\sum_{t=1}^n {}V_s(k)'V_t(k) \bigg  ) \nonumber \\
&=&\frac{1}{n} \sum_{s=1}^n \E_0\bigg (V_s(k)'V_s(k) \bigg  ) \nonumber \\
& = & \E_0\big(\|V_1(k)\|^2\big).
\en
We now remark that
\be*
\E_0\big (\|V_1(k)\|^2\big) & = & k\bigg(\di\frac{1}{k} \di\sum_{i=1}^k \E_0\big(Z_i^2 \big)\bigg)
\\
&\leq&
k\bigg(\di\frac{1}{k} \di\sum_{i=1}^k \E_0\big( Z_i^2 \big)^2\bigg)^{1/2}
\\
&\leq& k M^{1/2}.
\e*
Finally, 
we have
 $$\P_0(S_n\geq 2) \leq
 \sup_{1\leq k \leq d(n)}\bigg(\di\frac{1}{\E_0(\w{\lambda}_{\min}(k))^{1/2}}\bigg) \bigg(
\di\frac{{M^{1/4}}d(n)}{\sqrt{\log(n)}}\bigg).$$
 Theorem \ref{theo2} obtains as soon as we have shown that $\E_{0}(\w{\lambda}_{\min}(k))_{k>0}$ is a decreasing sequence,
 which is clear
 since matrices $(\E(\w{\S}_n(k))_{k>0}$ are embedded by construction, that is, the $k \times k$ submatrix   obtained from the $k$ first lines and $k$ first columns of $\w{\S}_n(k+1)$ coincides (in distribution) with $\w{\S}_n(k)$.

\cqfd

Finally, the test procedure is consistent against any
alternative having the form
\be*
H_1(q) : \exists q \in \N {\rm \ such \ that \ } a_i=b_i, \forall i=1,\cdots, q-1, {\rm \ and \ }  a_q \neq b_q,
\e*
where $a_i$ and $b_i$ are given by (\ref{basic2})
\begin{prop}
\label{prop2}
Under $H_1(q)$, ${T}_{n}(S_n)$ tends to infinity (in probability)
as $n \rightarrow \infty$.
\end{prop}
\textbf{Proof}
First note that $\lim_{n\rightarrow \infty}d(n) > q$.  We now prove that $\lim_{n\rightarrow \infty}P(S_n < q)= 0$.
For $k < q$ we have $P(S_n=k) \leq P(T_n(k) \geq T_n(q))$.
By the law of large numbers,
the variable $J_q/\sqrt{n}$ converges in probability to a non-null vector.
Since $\w{\Sigma}(q)$ is a positive definite matrix  we have
$J_q'\w{\Sigma}(q)^{-1}J_q={\cal O}_{\cal P}(n)$ and
the test statistics $T_n(q)-q\log(n)$ tends to $+\infty$ in probability
 under $H_1'$. By similar  arguments, $T_n(k) -k\log(n)$ tends to $-\infty$ under $H_1'$.
Then   $P(S_n=k)\rightarrow 0$ for all $k < q$.  It follows that $\lim_{n\rightarrow \infty}P(S_n\geq q)=1$ and then
$T_n({S_n})$ tends to $+\infty$ as $n\rightarrow  \infty$.

\cqfd

\begin{rmk}
It is well known  that the sample covariance matrix performs poorly in the high dimensional setting.
For applications in this context,  we could change the sample covariance  $
\w{\S}$ by a more suitable one.
In Ledoit and  Wolf (2004) a linear shrinkage is proposed, $\S^* = \rho_1 I + \rho_2 \w{\S}$,
where $I$ stands for the identity matrix and $\w{\S}$ for the sample covariance (like the one used in our paper).
Won et al. (2009) proposed a non linear shrinkage for Gaussian variance matrices.
Writing the sample covariance matrix $\w{\S} =Q diag(l_1,\cdots, l_p) Q^T$, their estimator has the form
$\S^* = Q diag(\w{\lambda_1},\cdots, \w{\lambda_p})Q^T$, where the $\w{\lambda}$'s are
constrained estimated  eigenvalues. 
Another approach is the thresholding procedure proposed in Cai and Liu
(2011,  see also El Karoui, 2008).
Writing $\w{\S}=(\w{\s}_{ij})_{k\times k}$
the sample covariance matrix, a universal thresholding estimator is $\S^*$ with
$\s_{ij}^*=\w{\s}_{ij}\mathbb I\{\w{\s}_{ij}\geq L_n\}$, with a proper choice of the threshold  $L_n$. Cai and Liu (2010) proposed the more adaptative thresholds
$L_{ij}=\delta (\theta_{ij} \log k /n)^{1/2} $, with  tuning parameter $\delta$ and some fourth  moment estimators $\theta$'s.  In our problem $\theta_{ij}$ should be the estimator of
the variance of $Var(Z_iZ_j)$.
\end{rmk}

\section{Numerical study}

\paragraph{Models and alternatives}

 We present empirical powers of the test through several models. We will denote by ${\cal P}(m)$ the Poisson distribution with mean $m$, ${\cal N}(a,b)$ the normal distribution with mean $a$ and standard error $b$, ${\cal B}(N,p)$ the binomial distribution on $[0,N]$ with probability of success $p$, $\chi^2_k$ the chi-squared distribution with degree $k$.
 We consider four models under $H_0$ and seven associated alternatives as follows:
 \bit
 \item
 Model MOD1:
 $Y \sim \chi^2_2$, $Z \sim {\cal N}(0,2)$, $V \sim \chi^2_2$, $W \sim {\cal N}(0,0.1)$. \\
 Alternative A11:
 $V \sim \chi^3_2$, $W \sim {\cal N}(0,0.1)$.
  \item
  Model MOD2:
 $Y \sim \chi^2_2$, $Z \sim {\cal N}(0,2)$, $V \sim \chi^2_2$, $W \sim {\cal N}(0,1)$. \\
  Alternative A12:
 $V \sim \chi^3_2$, $W \sim {\cal N}(0,1)$.
  \item
  Model MOD3:
 $Y \sim \chi^2_2$, $Z \sim {\cal N}(0,2)$, $V \sim \chi^2_2$, $W \sim {\cal N}(0,2)$. \\
  Alternative A13:
 $V \sim \chi^3_2$, $W \sim {\cal N}(0,2)$.
  \item
 Model MOD4:
  $Y \sim {\cal B} (10, 0.5)$, $Z \sim {\cal P}(2)$, $V \sim {\cal B}(10,0.5)$, $W \sim {\cal P}(1)$.
  \\
  Alternative A21:  $V \sim {\cal B}(10,0.4)$ and $W \sim {\cal P}(1)$,
  \\
  Alternative A22:  $V \sim {\cal B}(10,0.6)$ and $W \sim {\cal P}(1)$,
  \\
  Alternative A23:  $V \sim {\cal B}(9,0.5)$ and $W \sim {\cal P}(1)$,
  \\
  Alternative A24:  $V \sim {\cal B}(11,0.5)$ and $W \sim {\cal P}(1)$.
 \eit
For all models and alternatives we consider i.i.d.  data $(X_1,U_1), \cdots, (X_n,U_n)$ generated from two convolution models  satisfying (\ref{convol}). It is assumed that $Z$ and $W$ have known distribution.

\paragraph{Empirical levels}

We compute the test
statistic based on a sample size $n=30, 50, 100 $ and $200$ for a theoretical level $\alpha=5\%$.
The empirical level of the test is  defined as the percentage of rejection of the
null hypothesis over $10000$ replications of the test statistic under the null hypothesis.
We have fixed $d(n)=10$ arbitrarily large enough since the selected order does not exceed $4$ in all our simulations.

Empirical levels are reported in Table \ref{level}  for a fixed asymptotic level equal to 5\%. It can be seen that all values are
close to the asymptotic limit, also for small sample size.

\begin{table}[H]
	\caption{Empirical levels for MOD1, MOD2, MOD 3  and MOD4
 with sample sizes $30,50, 100, 200$}
\label{level}
	\begin{center}
	\begin{tabular}{ccccc}\hline
	Model & $n=30$ & $n=50$ & $n=100$ & $ n=200$  \\
	\hline
	MOD1 & 4.70 &5.07 & 4.92 & 4.90 \\
MOD2 & 4.51 & 4.98 &  4.66  & 5.01 \\
MOD3 & 4.38 & 4.72 &  4.84  & 4.94 \\
MOD4 & 4.80 & 4.93 &  4.80  & 4.63\\
	\hline
	\end{tabular}
	\end{center}
	\end{table}
\FloatBarrier

\begin{figure}[H]
		\includegraphics[scale=1.0]{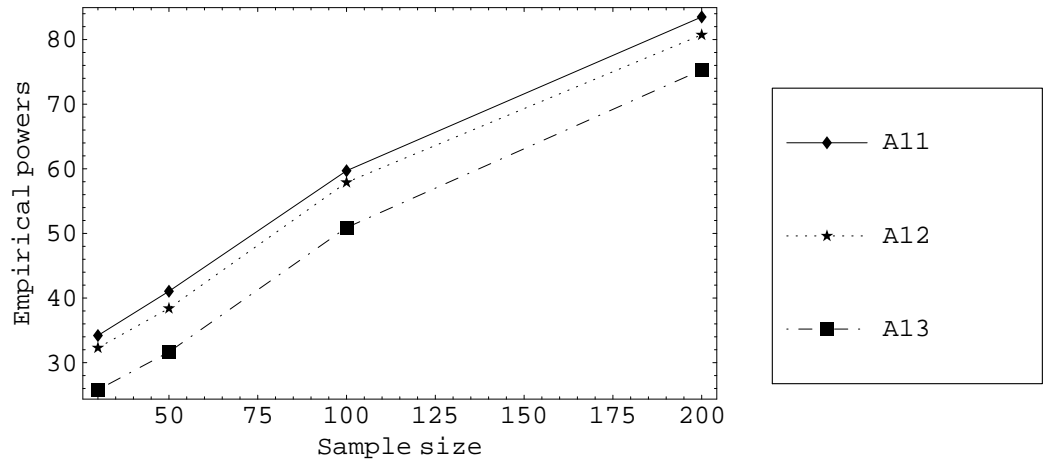}
		\caption{Empirical powers for alternatives A11 ($\blacklozenge$), A12 ($\bigstar$) and A13 ($\blacksquare$) with sample sizes $30,50, 100, 200$.
}\label{fig123}
	\end{figure}

\begin{figure}[H]
		\includegraphics[scale=0.8]{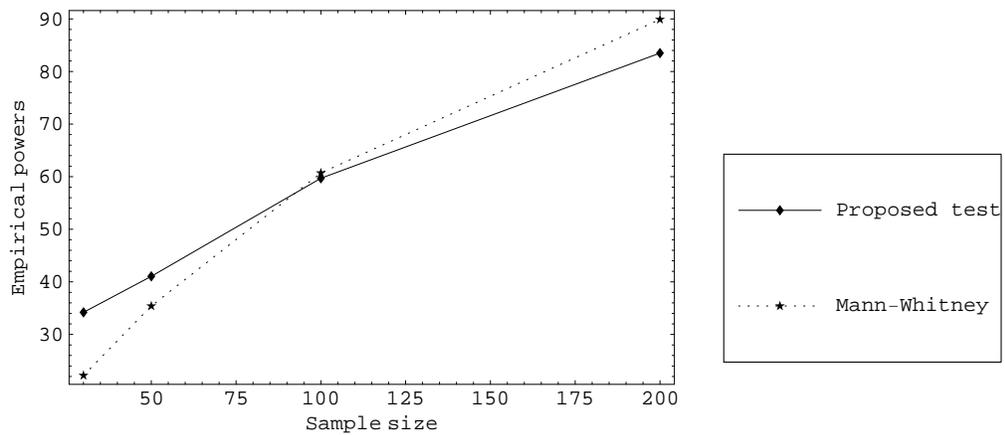}
		\caption{Empirical powers under alternatives A13 with the proposed method ($\blacklozenge$) and with the Mann-Withney
test ($\bigstar$) with sample sizes $30,50, 100, 200$.
}\label{fig2}
	\end{figure}

\begin{figure}[H]
		\includegraphics[scale=0.8]{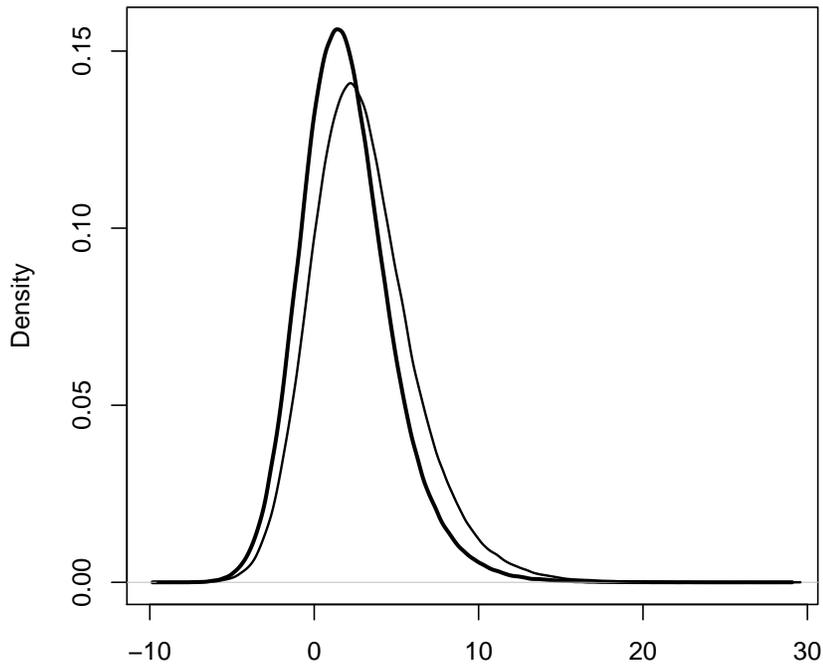}
		\caption{Density  under the null model MOD3 with $Y \sim \chi^2_2$ and $Z \sim {\cal N}(0,2)$ (left bold curve), and density under alternative A13 with $V \sim \chi^2_3$ and $W \sim {\cal N}(0,2)$ (right curve).
}\label{fig3}
	\end{figure}

\paragraph{Empirical powers}
The empirical power of the test is  defined as the percentage of rejection of the
null hypothesis over $10000$ replications of the test statistic under Alternative.
Empirical powers for alternatives A11-A13 are represented in Figure \ref{fig123}. In our knowledge, there is no equivalent
method in the literature to compare contaminated distributions and then it is not possible to confront these powers.
 However, for alternative A13 $Z$ and $W$ have the same distribution and the null hypothesis coincides with the equality of
 the two distributions ${\cal L}_X={\cal L}_U$. They consist in the convolution of a second order $\chi^2$ distribution with
 a Gaussian distribution ${\cal N}(0,2)$. Then, even if our method is not dedicated to the standard two-sample problem,
 we can compare its power with that of the classical Mann-Whitney test under A13.
 Figure \ref{fig2} shows that these two tests have similar powers, with a slight advantage for the Mann-Whitney test with
 large sample size. Note in Figure \ref{fig3} that the alternative A13 is close to a translation over the null distribution  that can be advantageous for the Mann-Whitney test.

\FloatBarrier

\FloatBarrier

\FloatBarrier

\begin{figure}[H]
		\includegraphics[scale=1.0]{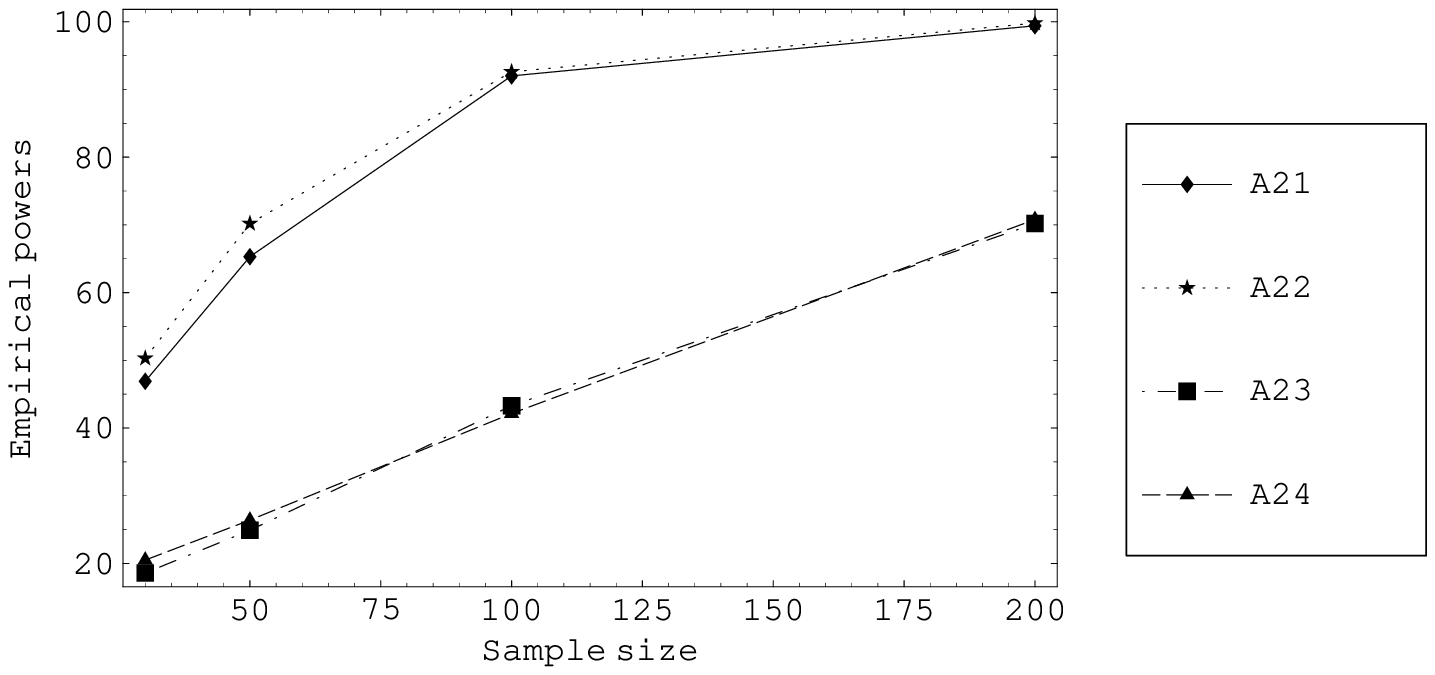}
		\caption{Empirical powers for alternatives A21 ($\blacklozenge$), A22 ($\bigstar$), A23 ($\blacksquare$) and  A24 ($\blacktriangle$) with sample sizes $30,50, 100, 200$.
}\label{fig4567}
	\end{figure}

\begin{figure}[H]
\vspace*{-1cm}
\subfigure[Distribution of $X$ under MOD4]
{
		\includegraphics[scale=0.7]{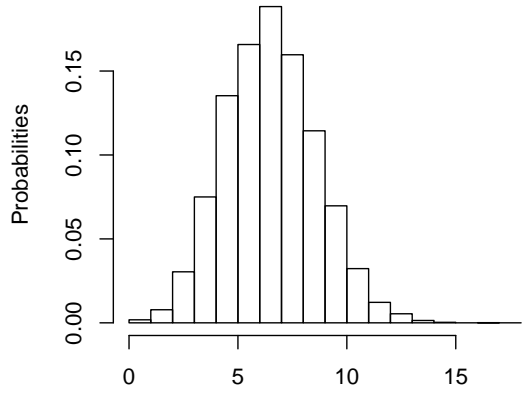}
}
\hspace*{-2cm}
\subfigure[Distribution of $U$ under alternative A22]
{
		\includegraphics[scale=0.68]{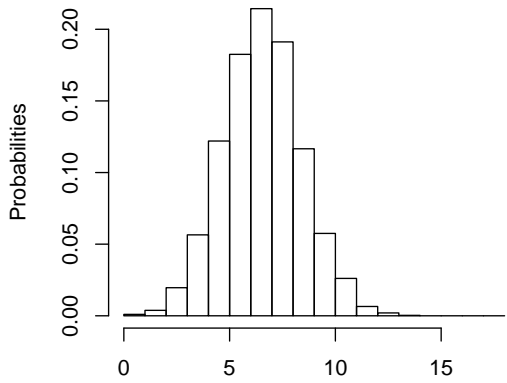}
}
\subfigure[Distribution of $U$ under alternative A23]
{		\includegraphics[scale=0.68]{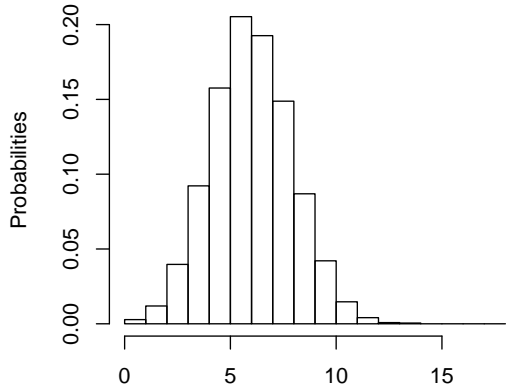}
}
\caption{Distributions under MOD4 (a), alternative A22 (b) and alternative A23 (c)}\label{fig4}
	\end{figure}

Figure \ref{fig4567} presents the powers of the test for MOD4 with alternatives A21-A24. Both alternatives A21 and A22 are very
well detected by the procedure. Under alternatives A23 and A24  the power is less good.  These results are essentially due to
the  nearness between the distributions of $Y$ and $V$ and not in that between $X$ and $U$. To illustrate this remark,
Figure \ref{fig4} shows  the proximity between $X$ and $U$ at once for alternative A22  and for alternative  A23.
All distributions are very similar.
But for alternative A22, the distributions of $Y$ and $V$ are closer than for alternative A23, explaining  its better  power.

\FloatBarrier

\FloatBarrier

\section{ Illustration}
Table \ref{uefa} reproduces paired data used by Meintanis (2007).
It concerns  matches of the UEFA Champion's League for the seasons 2004-05
and 2005-2006 where there was at least one goal scored by the home team, and there was at least one goal scored
directly from a kick by any team.
The first variable  $X$  is the time (in minutes) of the first kick
goal scored by either of the two team, and the second variable $U$ is the time of the first goal of any type scored by the home team.

\begin{table}[H]
\label{tabdata1}
\caption{UEFA Champion's league data}
\label{uefa}
~\\
\begin{tabular}{cccccc}
\hline
2005-2006 & $X$ & $U$ & 2004-2005 & $X$ & $U$ \\
\hline
Lyon-Real Madrid & 26    & 20 & Internazionale-Bremen & 34	& 34\\
Milan-Fenerbahce & 63	& 18 & Real Madrid-Roma &	53	& 39\\
Chelsea-Anderlecht & 19	& 19 & Man. United-Fernebahce &	54	& 7\\
Club Brugge-Juventus & 66	& 85 & Bayern-Ajax &	51	& 28\\
Fenerbhace-PSV & 40	& 40 & Moscow-PSG &	76	& 64\\
Internazionale-Rangers & 49	& 49 & Barcelona-Shakhtar &	64	& 15\\
Panathinaikos-Bremen &  8	& 8	 & Leverkusen-Roma &    26	& 48\\
Ajax-Arsenal & 69	& 71 & Arsenal-Panathinaikos  &	16	& 16\\
 Man. United-Benfica & 39	& 39 & Dynamo Kyiv-Real Madrid  &	44	& 13\\
Real Madrid-Rosenborg & 82	& 48 & Man. United-Sparta &	25	& 14\\
Villareal-Benfica  & 72	& 72 & Bayern-M. Tel-Aviv &	55	& 11\\
 Juventus-Bayern & 66	& 62 & Bremen-Internazionale &	49	& 49\\
 Club Brugge-Rapid & 25	& 9	 & Anderlecht-Valencia &    24	& 24\\
Olympiacos-Lyon & 41	& 3	 & Panathinaikos-PSV &    44	& 30\\
Internazionale-Porto & 16	& 75 & Arsenal-Rosenborg  &	42	& 3\\
Shalke-PSV & 18	& 18 & Liverpool-Olympiacos  &	27	& 47\\
Barcelona-Bremen & 22	& 14 & M. Tel-Aviv-Juventus&	28	& 28\\
Milan-Shalke & 42	& 42 & Bremen-Panathinaikos &	2	& 2 \\
Rapid-Juventus& 36	& 52 & &	&\\
 \hline	
\end{tabular}
\end{table}
\FloatBarrier

These data have been explored assuming a  Marshall-Olkin distribution in Meinatis (2007) and with a  bivariate generalized
exponential distribution in Kundu and Gupta (2009).  In Meintanis (2007) the conclusion was that   the Champion's-League
data may well have arisen from a  Marshall-Olkin distribution.  In Kundu and Gupta (2009) the generalized exponential
distribution can not be rejected for the marginals and the bivariate generalized exponential distribution can be used for
these data.
We consider here another model through contaminated Poisson distributions.

\paragraph{First model} First we assume an additive noise
\beg
\label{model1}
X = Y+Z & {\rm and }&  U=V+W,
\en
with $Y$ and $V$ having Poisson distributions and  $Z$ and $W$ being dependent random  noise with $\E(Z)=\E(W)=0$.
This model can be viewed as a mixed model with $Z$ and $W$ as paired random effects. These effects can be considered as discrete or continuous as in Meintanis (2007) or   Kundu and Gupta (2009).
 We assume that  $Y$ and $V$ have mean  (estimated) $40.9$ and $32.9$. The observed variances are  larger than the means thereby believe there is a phenomenon of overdispersion.
 Obviously under (\ref{model1}) we have
 $Var(X) = Var(Y) +Var(Z)$ and $Var(U)= Var(V)+Var(W)$.
 We apply our procedure to test the equality of the distributions of $Z$ and $W$.
 \\
 {\it Conclusion:} The first statistic $T(1)$ is retained and
  we obtain a p-value equal to $0.28$.  Hence there is no evidence that the two additive paired random effects differ.

\paragraph{Second model}
We also consider a multiplicative noise yielding to the following scale mixture
\beg
\label{model2}
X = YZ & {\rm and }&  U=VW,
\en
with $Y$ and $V$ having Poisson distributions and  $Z$ and $W$ being real positive dependent random scale factor with $\E(Z)=\E(W)=1$.
Again this model  can be viewed as a mixed model with random paired effects. The observed values are discretized but we can assume that $Z$ and $W$ are discrete or continuous.
 We assume that  $Y$ and $V$ have mean  (estimated) $40.9$ and $32.9$ and there is still a phenomenon of overdispersion assuming it is a standard Poisson model.
 Under (\ref{model2}) the variances satisfy
 $Var(X) = (2Var(Z)+1)\E(X)$ and  $Var(U) = (2Var(W)+1)\E(U)$.
 Our purpose is to test $H_0: {\cal L}_Z = {\cal L}_W$, or equivalently $ {\cal L}_{\log(Z)} = {\cal L}_{\log(W)}$.
 For that we consider the transformation of (\ref{model2})
 \be*
 \log(Z) = \log(Y) + \log(Z) & {\rm and } & \log(U) = \log(V)+\log(W).
 \e*
{\it Conclusion:}  Using our method we obtain a p-value equal to $0.70$. Again we see that the multiplicative paired random effects seem to
 have the same distribution.

\section{Discussion}

This paper discusses the problem of comparing two distributions contaminated by different noises. The test is very simple
and allows to compare two independent as well as two paired contaminated samples.
Simulation studies suggest that the proposed method works well with an empirical level close to that expected.

It may be  noted that the test statistic is decomposed into moments of $X,Z,U,$ and $W$.
Then it is clear that only the knowledge of
the moments of $Z$ an $W$ are required  instead of their  distributions.
Hence  the test could be adapted when these distributions are unknown, if    their moments can be estimated
from independent samples.

Eventually, the multivariate case could be envisaged by using the following characteristic property:
if $Y$ and $V$ are two random vectors taking values in $\R^d$ then we have
\be*
H_0 : Y =^d V & \Leftrightarrow & \forall \|u\|\leq 1, u'Y =^d u'V,
\e*
and clearly multidimensional observations can be transformed into unidimensional ones by applying a sequence of vectors $u$ on $X$ and $U$.
For a fixed value of $u$ the problem consists in an univariate test and the statistic
$T_n(S_n)$ can be used. denoting by $T_n(u)$ this statistic the
process $\{T_n(u) ; u \in (0,1)^d \}$ converges to a  Gaussian process and a new test statistic can be envisaged
by estimating the covariance  operator of the process to get a  $\chi^2$  null distribution.
In practice the sequences of vectors $u$ can be randomly chosen in $(0,1)^d$, but it can also be done by a Quasi Monte
Carlo method (see for instance L'Ecuyer, 2006).

To conclude, the multisample case can also be envisaged as follows:
assume that we have $d$ convolutions simultaneously
\be*
X(i)  =  Y(i) +  Z(i), & & i=1, \cdots, d
\e*
observed from $d$ samples.
Write $\alpha_j(i)=\E(Y(i)^j)$ and $\bar{\alpha}_j = \frac{1}{d}\sum_{i=1}^d \alpha_j(i)$ the common
value under the null hypothesis $H_0: {\cal L}_{Y(1)}= \cdots = {\cal L}_{Y(d)}$.
Then under $H_0$ the $k\times d$ vector  $D$ with components $D_{ij}=\alpha_j(i)-\bar{\alpha}_j$ is centered and normally
distributed. An adaptation of the data driven smooth test seems then possible.


\end{document}